\documentclass[12pt]{article}

\usepackage{amsfonts}
\usepackage{amsmath,amssymb}
\usepackage{amscd}
\usepackage{amsthm}
\newtheorem{df}{Definition}[section]

\newtheorem{lem}{Lemma}[section]
\newtheorem{cly}{Corollary}[section]
\newtheorem{prop}{Proposition}[section]
\newtheorem{rem}{Remark}[section]
\newtheorem{ex}{Example}[section]

\title{On unitality conditions for hom-associative algebras}
\author {Ya\"el Fr\'egier; Aron Gohr\\{\small Mathematics Research Unit}\\
{\small  162A, avenue de la faiencerie}\\
{\small  L-1511, Luxembourg}\\
 {\small  Grand-Duchy of Luxembourg.}\\
}
\begin {document}
\date{}
\maketitle

\begin{abstract}
In \emph{hom-associative} structures, the associativity condition $(xy)z=x(yz)$ is twisted to $\alpha(x)(yz) = (xy)\alpha(z)$, with $\alpha$ a map in the appropriate category. In the present paper, we consider two different unitality conditions for hom-associative algebras. The first one, existence of a unit in the classical sense, is stronger than the second one, which we call \emph{weak unitality}. We show associativity conditions connected to the size of the image of the twisting map for unital hom-associative algebras. Also the problem of embedding arbitrary hom-associative algebras into unital or weakly unital ones is investigated. Finally, we show that weakly unital hom-associative algebras with bijective twisting map are twisted versions of associative algebras. 
\end{abstract} 


\section* { Introduction.}
We fix first some conventions and notations. In this article, $k$ will be a commutative ring, $K$ a field. Modules and algebras will be understood to be over an arbitrary commutative ring. If $\alpha:G \rightarrow H$ is a homomorphism of groups (rings, modules, etc.) we will denote by $Ke(\alpha)$ its kernel and by $Im(\alpha)$ its image. $V$ will be a $k$-module.\\

Hom-algebras were first introduced, in the Lie case, by Hartwig, Larson, Silvestrov \cite{HLS}. This notion was later extended by dropping conditions on the twisting and transferred to the associative category by Makhlouf, Silvestrov \cite{MSI}. Hom-algebraic counterparts have been found for many classical algebraic constructions, see e.g. \cite{AMS}, \cite{LS1}, \cite{MSII}, \cite{MSIII}, \cite{MSIV}, \cite{Yau:EnvLieAlg}, \cite{Yau:HOMology} and \cite{Yau:Hom-bi}. The present work is aimed towards supporting these efforts by contributing to a better understanding of the structure theory of hom-associative algebras.\\

This paper is divided into two sections. In Section \ref{section1-uha} we consider \emph{unital} hom-associative algebras. More specifically, in (\ref{elementary properties}) we recall the definitions of hom-associative rings and algebras and derive some elementary but for our studies very useful results. The most important of these is Prop. \ref{imassociative}, which gives various association conditions for elements in the image of the twisting map of a unital hom-associative algebra. It is essentially the observation contained in Prop. \ref{imassociative} which motivated our study of associativity conditions for unital hom-algebras. The results obtained in Section \ref{elementary properties} will also be useful in a subsequent paper, where we investigate certain variations on the theme of hom-associative algebras.\\
In (\ref{associativity conditions}), we derive first some corollaries of Prop. \ref{imassociative}. We then turn to the special case of unital hom-associative algebras over fields in Prop. \ref{codim}. The assumption of the base ring being a field is important here because the associativity criteria we derive involve the codimension of $Im(\alpha)$, viewed as a vector subspace of our hom-algebra $A$, where $\alpha$ is the twisting map of $A$. We close subsection \ref{associativity conditions} by discussing some open questions and counterexamples to associativity conditions that do not hold.\\
In (\ref{twisting structure}) we apply the results obtained previously to the following in our view natural question:\\

\emph{Suppose $(A, \star)$ is a not necessarily associative, unital algebra. Under which twisting maps does this become a hom-associative algebra?}\\

We give an answer to this question in Prop. \ref{structtwist}. Essentially, we identify a one-one correspondence between such twistings and elements of the center of our algebra which satisfy some additional association conditions.\\

In Section \ref{weak}, we relax the condition of unitality. Our investigation is motivated by the fact that there are, relatively speaking, much less unital specimens in the class of hom-associative algebras than there are unital associative algebras among the associative ones. We also observe that in some sense, the unit element plays a double role in a unital associative algebra when it is viewed as hom-associative. A decoupling of these two roles of the unit leads, in subsection \ref{weak} to the introduction of \emph{weakly} unital hom-associative algebras. We answer in the negative the question of whether any hom-associative algebra can be embedded into at least a weakly unital one, and point out, in Ex. \ref{yau-example}, that all hom-associative algebras obtained by a procedure due to Yau \cite{Yau:HOMology} are weakly unital.\\

The rest of the paper is devoted to a study of the relationship between Ex. \ref{yau-example} and weakly unital hom-associative algebras. Our main finding, Prop. \ref{weakly unital algebras are twistings}, says that \emph{all} weakly unital hom-associative algebras with bijective twisting map can be viewed as twistings of unital associative algebras by a generalization of Yau's procedure. The idea of the proof is that, in the case of a hom-associative algebra $(A, \star, \alpha)$ with a bijective twisting map $\alpha$, the twisting procedure in Yau \cite{Yau:HOMology} can be applied to $A$ itself with twisting map $\beta := \alpha^{-1}$ to obtain a new $k$-algebra $(A, \star')$. By construction, if the twisting procedure is applied to $(A, \star')$ with twisting map $\alpha$, one recovers the original algebra $(A, \star, \alpha)$. Surprisingly, if $(A, \star, \alpha, c)$ was weakly unital, $(A, \star')$ is always associative and unital.\\



 \section{Unital hom-associative algebras} \label{section1-uha}
\subsection{Elementary properties} \label{elementary properties}
The object of study of this section are \emph{unital} hom-associative rings and algebras. For precision, we set the following definition, following \cite{MSI}:
\begin{df} \label{homring-definition}
Let $(V, +, \star, \alpha)$ be a set together with two binary operations $+$ and $\star$ as well as one unary operation $\alpha:V \rightarrow V$. Then $(V, +, \star, \alpha)$ is called a \emph{hom-ring} if $(V, +, \star)$ is a not necessarily associative ring and $\alpha$ is an abelian group endomorphism of $(V,+)$ such that the hom-associativity condition
\begin{displaymath}
\alpha(x) \star (y \star z) = (x \star y) \star \alpha(z)
\end{displaymath}
is fulfilled for all $x,y,z \in V$. $V$ is called \emph{unital} if there exists an element $1 \in V$ such that $1 \star x = x \star 1 = x$ for all $x \in V$. Left-unital and right-unital hom-rings and hom-algebras are defined accordingly.
\end{df}

We start our investigations by a series of lemmas which are useful for carrying out calculations in the unital hom-associative setting:



\begin{lem}\label{lemma1}
Let $(V,\star,\alpha,1)$ be a unital hom-associative algebra. One has for all
$x,y$ in $V$:
\begin{eqnarray}
\alpha(x)\star y & = & x\star \alpha (y) \label{eq-aa}\\
x\star \alpha (1) & = & \alpha (x), \label{eq-ee} \\
\alpha (x\star y) & = & x\star \alpha (y). \label{eq-ii}
\end{eqnarray}
\end{lem}

\begin{proof} Left as an easy exercise to the reader.
\end{proof}
Some other remarks on (Lemma \ref{lemma1}) are in order. First, analogously to Eq. (\ref{eq-ee}), one can also show that $\alpha(x) = \alpha(1) \star x$, so $\alpha(1)$ commutes with all elements. On the other hand, if $(V, \star)$ is an associative algebra, multiplication with an arbitrary central element induces a hom-associative structure. These two observations taken together allow a complete classification of hom-associative structures compatible with $(V, \star)$ if $(V, \star)$ is a unital associative algebra: the set of admissible twisting homomorphisms $\alpha$ corresponds one-to-one with the center of $V$.\\
 Second, Eq. (\ref{eq-ii}) shows that $Im(\alpha)$ and $Ke(\alpha)$ are stable under multiplication with arbitrary elements. Since they are also clearly stable under $\alpha$, this means that $Im(\alpha)$ and $Ke(\alpha)$ are hom-algebra ideals. In fact, it is an immediate consequence of the following proposition that $(Im(\alpha), \star)$ is an associative algebra:




\begin{prop}\label{imassociative} Let $(V, \star, \alpha, 1)$ be a unital hom-associative algebra. One has for all $x,y$ and $z$ in $V$:
\begin{eqnarray}
\alpha(x)\star (y\star z) & = & (\alpha(x)\star y)\star z, \label{eq-1}\\
x\star (\alpha(y)\star z) & = & (x\star \alpha(y))\star z, \label{eq-2}\\
x\star (y\star \alpha(z)) & = & (x\star y)\star \alpha(z), \label{eq-3}\\
\alpha(x \star (y \star z)) & = & \alpha((x \star y) \star z). \label{eq-4}
\end{eqnarray}
\end{prop}

\begin{proof} The proof is an easy application of the results of lemma \ref{lemma1} and the definition of a unital Hom-associative algebra and is also left as an exercise.
 
\end{proof}

\begin{rem}
One can rephrase equations Eq. (\ref{eq-1}-\ref{eq-3}) by saying that the image of $\alpha$ lies in the \it{nucleus} of the non-associative algebra $(V, \star)$, i.e. in the set of elements that associate with all elements of $V$.
\end{rem}

\subsection{Associativity conditions for unital hom-associative algebras} \label{associativity conditions}

A natural question is the following: "When is a hom-associative
algebra really associative?". A first answer is a trivial corollary
of the previous proposition:
\begin{cly}\label{surjectivity implies associativityI}
$(V, \star)$ is associative if $\alpha$ is surjective.
\end{cly}

 If as base ring we take the field $K$ and if we restrict ourselves to the finite dimensional case, $\alpha$ is injective exactly if surjective by linear algebra. Hence the preceding remark shows then that $(V, \star)$ is associative whenever $\alpha$ is injective. We will see in a minute that, in fact, \emph{surjectivity of $\alpha$ implies injectivity} in general and that \emph{injectivity of $\alpha$ forces associativity}:
\begin{cly} \label{injectivity implies associativity}
$(V, \star)$ is associative if $\alpha$ is injective.
\end{cly}
\begin{proof} For all $x,y,z \in V$ we have $\alpha((x \star y) \star z) = \alpha(x \star (y \star z))$ by the previous theorem. But if $\alpha$ is injective, this means associativity.
\end{proof}
Now assume $\alpha$ is surjective. Then there is some $a \in V$ such that $\alpha(a) = a \alpha(1) = 1$. Recall that $\alpha(1)$ is in the center of $V$. Then if $\alpha(x) = 0$ for some $x \in V$, we see $0 = a  \star \alpha(x) = a \star (x \star \alpha(1)) = a \star (\alpha(1) \star x) \overset{\text{Eq. (\ref{eq-2})}}{=} x$. We have hence proven:
\begin{cly} \label{surjectivity implies injectivity}
If $\alpha$ is surjective, then $\alpha$ is also injective.
\end{cly}
The inverse implication is wrong, because e.g. with $A$ a commutative (and associative) algebra and $r \in A$ a non-zerodivisor one can obtain with $\alpha(x) = rx$ an injective $\alpha$ which induces on $A$ a hom-associative structure. But this $\alpha$ will not in general be surjective.\\
So we know that surjective $\alpha$ implies associativity of $(V, \star)$ and that in fact the weaker condition of injectivity is sufficient. It is natural then to ask whether we might generalize the observation that a surjective $\alpha$ induces associativity in other ways. For instance, we can replace surjectivity by the following weaker conditions on codimension:
\begin{prop}\label{codim} Let $(V, \star, \alpha)$ be a unital hom-associative algebra over a field $K$. If any of the following conditions is satisfied, $(V, \star)$ is associative:
\begin{enumerate}
\item $codim(Im(\alpha)) \leq 1$
\item $codim(Im(\alpha)) \leq 2$ and $(V, \star)$ commutative
\item $codim(Im(\alpha)) \leq 2$ and $\alpha$ injective on $Im(\alpha)$.
\end{enumerate}
\end{prop}
\begin{proof}
1. Suppose $(V, \star)$ is not associative. Then $\alpha$ is not surjective by Corollary \ref{surjectivity implies associativityI}, hence $codim(Im(\alpha)) = 1$ and $1 \not\in Im(\alpha)$. But this means with the usual embedding $K \rightarrow V$ that $V = K \oplus Im(\alpha)$. The expressions $x_1 \star (x_2 \star x_3)$ and $(x_1 \star x_2) \star x_3$ coincide for any $x_1, x_2, x_3 \in V$ such that for at least one $i$ we have $x_i \in K$ or $x_i \in Im(\alpha)$ due to Proposition \ref{imassociative}. Association of arbitrary elements of $V$ follows from this. Note that without the condition on codimension, this argument still shows that $K \oplus Im(\alpha)$ is a hom-subalgebra of $V$ which contains only elements that associate with the rest of $V$ and among each other.\\
2. By the previous argument, we know that $codim(Im(\alpha)) = 2$ if $(V, \star)$ is non-associative and $codim(\alpha) \leq 2$. Choose a direct complement $U \subseteq V$ of $K \oplus Im(\alpha)$ in $V$, then by $codim(Im(\alpha)) = 2$ we see that $V = K \oplus Im(\alpha) \oplus U$ and that $dim_K(U) = 1$. Since products $x \star y \star z$ with at least one of $x,y,z$ an element of $V' = K \oplus Im(\alpha)$ are associative, the only thing that needs to be checked is then $u \star (u \star u) = (u \star u) \star u$, where $u$ is a $K$-linear generator of $U$. But this is obvious given commutativity.\\
3. We analyze now the remaining case, where $\alpha$ is assumed injective on $Im(\alpha)$. Using the same notations as in (2), we see that again we have to prove $u \star (u \star u) = (u \star u) \star u$.
To see this, we introduce for any $x \in V$ the notation $x =: (x_K, x_\alpha, \lambda_x)$ with $x_K \in K, x_\alpha \in Im(\alpha), \lambda_x \in K$ defined by the decomposition $x = x_K + x_\alpha + \lambda_x u$. Then, we get $u \star u = (\lambda_K, \lambda_\alpha, \lambda)$ and therefore
\begin{displaymath}
u \star (u \star u) = (\lambda_K \lambda, u \lambda_\alpha + \lambda \lambda_\alpha, \lambda_K + \lambda^2)
\end{displaymath}
and
\begin{displaymath}
(u \star u) \star u = (\lambda \lambda_K, \lambda_\alpha u + \lambda \lambda_\alpha, \lambda_K + \lambda^2).
\end{displaymath}
One uses in these calculations that $Im(\alpha)$ is an ideal and that vector subspaces of $V$ are of course stable under multiplication with elements of the base field. We see then that $u \star (u \star u) - (u \star u) \star u = u \lambda_\alpha - \lambda_\alpha u \in Im(\alpha)$. Now $\alpha(u \star (u \star u) - (u \star u) \star u) = 0$ because of Eq. (\ref{eq-4}), so by injectivity of $\alpha$ on $Im(\alpha)$ the desired equation $u \star (u \star u) = (u \star u) \star u$ is obtained.
\end{proof}
It is natural to ask if $codim(Im(\alpha)) = 2$ is sufficient in general to force associativity. We do not know a proof of this nor a counterexample, although the last calculation of the previous proof in our view suggests that this is not generally so.  It is clear that in codimension three or higher, hom-associativity does no longer imply associativity, since there exist unital non-associative algebras of vector space dimension three, and with the zero map as twisting they become hom-associative.\\
In finite dimension, the condition $codim(Im(\alpha)) = 2$ is equivalent to\\ $dim(Ke(\alpha)) = 2$ by linear algebra. While we do not know whether the former condition is enough to force associativity, we can show that the latter is not:
\begin{ex} \label{counterex dim two conjecture}
Let $K$ be a field and let $A := K[X]$. Suppose further that $U$ is a two-dimensional $K$-vector space on which a non-associative product is defined. We can make $U$ into an $A$-module through the ring homomorphism $K[X] \rightarrow K$ induced by evaluation at zero. We set $B := A \times U$ and define on $B$ the multiplication
\begin{displaymath}
(a_1, u_1) \star (a_2, u_2) := (a_1 a_2, a_1 u_2 + a_2 u_1 + u_1 u_2).
\end{displaymath}
We set further $\alpha(a,u) := (Xa, Xu) = (Xa, 0)$. Then straightforward calculation shows that $Ke(\alpha) = U$ and that $(B, \alpha)$ is hom-associative with unit $(1,0)$. But $B$ is also non-associative because $U$ was. So $dim(Ke(\alpha)) = 2$ does not force associativity.
\end{ex}
Unfortunately, there is no analog of this example in finite dimension, where it could be used to conclude that also $codim(Im(\alpha)) = 2$ does not imply associativity. The reason is that with $A$ finite dimensional as a vector space and $a$ an element in the center of $A$ annihilating $U$, $a$ would either be a zero divisor or an invertible element, since $A$ is artinian if finite dimensional. In the first case, $U$ is a proper subset of the kernel of $\cdot a$ in $B$. The second case contradicts the assumption that $a$ annihilates $U$.\\
We do not know if $dim(Ke(\alpha)) = 1$ forces associativity in the infinite-dimensional case.\\
We will quickly discuss now another interpretation of our findings around these topics. We have seen that $Ke(\alpha)$ and $Im(\alpha)$ are hom-ideals of the hom-associative algebra $(V, \star, \alpha, 1)$. Proposition \ref{imassociative} also tells us that for arbitrary $x,y,z \in V$, we have $x \star (y \star z) - (x \star y) \star z \in Ke(\alpha)$. But this means that the hom-algebra $A := V/Ke(\alpha)$ is associative and, more specifically, hom-associative under the injective twisting map induced on $A$ by $\alpha$. We call $A$ the \emph{associative factor} of $V$. If $\alpha$ is nonzero, it is unital. By definition, the operations on $A$ recover, in a sense, the operations on $V$ up to a perturbation in $Ke(\alpha)$.\\

\subsection{Structure of twisting maps.} \label{twisting structure}
Here is now an application of the results shown previously. It answers the
question: "which hom-associative structures may be installed on a
given algebra?".
\begin{prop}\label{structtwist} Suppose that $k$ is a commutative ring and that $(A, \star, 1)$ is a unital $k$-algebra. Then the twisting maps $\alpha:A \rightarrow A$ which are compatible to $A$ in the sense that $(A, \star, \alpha, 1)$ becomes a hom-associative algebra form a commutative associative algebra with respect to pointwise addition of functions as additive operation and composition as multiplication. They correspond one-to-one with elements $a \in A$ such that:
\begin{enumerate}
\item $ax = xa$ for all $x \in A$ and
\item the set $Aa$ forms an ideal in $A$ and every element of this ideal associates with every element of $A$. 
\end{enumerate}
Denoting the set of compatible twisting maps by $Twist(A)$ and the set of elements in $A$ which satisfy the preceding conditions by $AC(A)$, the correspondence is given by $\Phi:Twist(A) \rightarrow AC(A)$ and $\Psi: AC(A) \rightarrow Twist(A)$, defined through 
\begin{eqnarray*}
\Phi(\alpha) := \alpha(1),\\
\Psi(x)(y) := x \star y.
\end{eqnarray*}
\end{prop}
\begin{proof}
We have already proven that in a unital hom-associative algebra $(V, \star, \alpha, 1)$, the element $\alpha(1)$ has all of the desired properties. On the other hand, assume that $a \in A$ satisfies the conditions (1,2). Then we have for arbitrary $x,y,z \in A$ that
\begin{displaymath}
(ax)(yz) \overset{1.}{=} (xa)(yz) \overset{2.}{=} x(a(yz)) \overset{1.}{=} x((yz)a) \overset{2.}{=} x(y(za)) \overset{2.}{=} (xy)(za)
\end{displaymath}
and therefore hom-associativity of $A$ with $\alpha(x) = ax$. By the definitions, we have
\begin{displaymath}
\Psi(\Phi(\alpha))(x) = \alpha(1) \star x = \alpha(x)
\end{displaymath}
and
\begin{displaymath}
\Phi(\Psi(x)) = x \star 1 = x
\end{displaymath}
for all $x$, so $\Psi$ and $\Phi$ are inverse to each other. It is clear that both maps are additive and one only needs to verify that they are also multiplicative. We calculate:
\begin{displaymath}
\Phi(\alpha_1 \circ \alpha_2) = \alpha_1(\alpha_2(1)) = \alpha_1(\Phi(\alpha_2)) = \Phi(\alpha_1) \Phi(\alpha_2). 
\end{displaymath}
Since a bijective additive map of (not necessarily unital) commutative rings which is multiplicative is a ring isomorphism, and since $AC(A)$ is clearly closed under multiplication, this concludes the proof.
\end{proof}






\section{Weakly unital hom-associative algebras} \label{weak}
\subsection{Embeddings of non-unital hom-algebras.} \label{embedding}
We start this section by studying the concept of unitalization of hom-associative algebras.
For $A$ an associative algebra, not necessarily unital, over a commutative ring $k$, it is not difficult to find an embedding $A \rightarrow A'$ such that $A'$ is a unital $k$-algebra. One can simply set $A' := k \oplus A$ and use the multiplication rule $(\lambda, a) \cdot (\mu, b) := (\lambda \mu, \lambda b + \mu a)$. The existence of such an embedding implies in particular that identities which can be proven to hold for tuples of arbitrary elements of arbitrary \emph{unital} associative algebras must also hold for not necessarily unital associative algebras.\\
In the case of hom-associative algebras, this is not so. For instance, the following example can be used to show that in a non-unital hom-associative algebra, the identity $\alpha(x) \star y = x \star \alpha(y)$ does not necessarily hold:
\begin{ex} \label{non_adjoint} Let $K$ be a field and set $A = K^2$ with multiplication
\begin{displaymath}
(\lambda_1, \lambda_2) \star (\mu_1, \mu_2) := (0, \lambda_1 \mu_1).
\end{displaymath}
Define further $\alpha:K^2 \rightarrow K^2$ by
\begin{displaymath}
\alpha(\lambda, \mu) := (\lambda + \mu, \mu).
\end{displaymath}
Then for all $x,y,z \in V$ we find $\alpha(x) \star (y \star z) = (x \star y) \star \alpha(z) = 0$, so $(A, \star, \alpha)$ is hom-associative. In fact, $(A, \star)$ is even a nonunital associative commutative algebra. But, $\alpha(0,1) \star (1,0) = (0,1)$ is not equal to $(0,1) \star \alpha(1,0) = (0,0)$.
\end{ex}
One can infer from failure of such identities to hold in the nonunital case that there cannot exist an unitalization procedure for hom-associative algebras which works in every case 
In the last example, one notes that a unit can easily be adjoined to $A$ if it is endowed with a different hom-associative structure, since we can view it e.g. as an associative algebra or make it trivially hom-associative by setting $\alpha = 0$. On closer examination, we realize that the unit element of a unital associative algebra plays two different roles if we view the associativity condition as a particular instance of hom-associativity. On the one hand, the unit is simply the neutral element of multiplication; but on the other hand, in a unital hom-algebra, the twisting homomorphism $\alpha$ is necessarily given by multiplication with some element from the algebra itself, and this is the unit element in case of $\alpha = id$.\\
This motivates the following definition:
\begin{df}
Let $(A, \star, \alpha)$ be a hom-associative algebra. Then $A$ is called \emph{left weakly unital} if $\alpha(x) = cx$ for some $c \in A$. Analogously, it is called \emph{right weakly unital} if $\alpha(x) = xc$ for some $c \in A$. $A$ is called \emph{weakly unital} if it is both left and right weakly unital. In these cases, the element $c \in A$ is called a \emph{weak (left/right) unit} of $A$.
\end{df}
This notion has been introduced in a preprint version of this paper. One can remark that since then it has already been used in several papers by other authors e.g. \cite{Yau:Hom-Yang}, \cite{Yau:Hom-qgI}, \cite{Yau:Hom-qgII} and \cite{CG09}.
\begin{rem}
In general, a hom-algebra can have many weak units if it has any. However, if $\alpha$ is injective and if $c_l$ is a weak left unit of $A$ and $c_r$ is a weak right unit of $A$, the two coincide, as in this case $\alpha(c_l) = c_l \star c_r = \alpha(c_r)$.
\end{rem}
From the theory we have built up so far, it is clear that unital hom-associative algebras are always weakly unital. However, the converse is evidently not the case, as for instance an arbitrary algebra if equipped with the trivial hom-structure is also weakly unital. More interesting examples of weakly unital but not unital hom-algebras can e.g. be obtained from the following construction due to Yau \cite{Yau:HOMology}:
\begin{ex} \label{yau-example}
Let $(A, \mu)$ be a unital associative algebra with multiplication $\mu$ and let $\alpha: A \rightarrow A$ be an algebra endomorphism. Set then for $x,y \in A$
\begin{displaymath}
x \star y := \alpha(\mu(x,y)).
\end{displaymath}
Then $(A, \star, \alpha)$ is a weakly unital, hom-associative algebra with weak unit $1$. 
\end{ex}
\begin{proof} Hom-associativity is proven in Yau \cite{Yau:HOMology}. Weak unitality is $1 \star x = \alpha(\mu(1,x)) = \alpha(x)$ and $x \star 1 = \alpha(x)$ in the same way. This concludes the proof.
\end{proof}
  One might hope, then, that every hom-associative algebra may be embedded into one that is weakly unital. The first main point we wish to make in this section is that this is not the case:
\begin{rem} There is no embedding of the algebra in Example \ref{non_adjoint} into a weakly unital hom-associative algebra.
\end{rem}
\begin{proof} Assume that $A \subseteq B$ with $(B, \beta)$ weakly unital, i.e. $\beta(x) = cx$ for some $c \in B$ and $\beta|_A = \alpha$. Denoting by $e_1, e_2$ the standard basis vectors of $K^2 = A$, we would then have $e_1 e_2 = 0$ and $\alpha(e_1) \alpha(e_2) = e_1 \neq 0$. But also by hom-associativity
\begin{displaymath}
\alpha(e_1) \alpha(e_2) = (c e_1) (\alpha e_2) \overset{\text{Hom}}{=} (c c) (e_1 e_2) = 0,
\end{displaymath}
so a contradiction is obtained.
\end{proof}
\subsection{Weakly unital hom-algebras with bijective twisting}
It is clear that among weakly unital hom-associative algebras, a case that is of particular interest is the one of injective $\alpha$. A number of special properties follow in this case, for instance as was already pointed out, weak units are uniquely defined. Also, the associativity constraint given by the hom-associativity condition is particularly strong when $\alpha$ is injective. Natural examples of such algebras are given for instance by the procedure given in Ex. \ref{yau-example} when the algebra endomorphism $\alpha$ is injective.\\
The rest of this section is devoted to the classification of weakly left unital hom-associative algebras with \emph{bijective} twisting map. The main result will be that all of them can be obtained from the following generalization of Ex. \ref{yau-example}:
\begin{ex} \label{generalized yau-example}
Let $(A, \cdot)$ be a $k$-algebra, not necessarily associative, and let $\alpha:A \rightarrow A$ be a $k$-linear map satisfying the equation
\begin{equation} \label{inducing homass}
\alpha(\alpha(x) \alpha(y z)) = \alpha(\alpha(xy) \alpha(z)).
\end{equation}
Set then
\begin{displaymath}
x \star y := \alpha(xy)
\end{displaymath}
for any $x, y \in A$. Then $(A, \star, \alpha)$ is a hom-associative $k$-algebra. If $1 \in A$ is a (left/right) unit of $A$, then $(A, \star, \alpha, 1)$ becomes weakly (left/right) unital.
\end{ex}
\begin{proof}
To check hom-associativity, we calculate
\begin{displaymath}
\alpha(x) \star (y \star z) = \alpha(\alpha(x) \alpha(yz)) = \alpha(\alpha(xy) \alpha(z)) = (x \star y) \star \alpha(z).
\end{displaymath}
If $1 \in A$ is a left unit, we see as in Ex. \ref{yau-example} that $\alpha(x) = 1 \star x$ for all $x \in A$. This concludes the proof.
\end{proof}

Our target for the rest of the paper will be the proof of the following statement:\\

\emph{All weakly unital hom-associative algebras with bijective twisting arise through the preceding example from associative algebras.\\}

 In fact we are going to show that a weak one-sided unit is sufficient. A recurring theme in the proof will be to show that certain algebraic structures satisfy conditions which may be viewed as variations of hom-associativity. This observation can be used to motivate a systematic study of hom-associativity conditions different from the one introduced in Makhlouf, Silvestrov \cite{MSI}, which we carry out in another paper. Eq. (\ref{inducing homass}) may be seen as a first example of such an alternative hom-associativity condition. Particularly striking examples are also the identities proven in Lemma \ref{induced type I3 lemma} and  Lemma \ref{weak unit main tool}.\\
For brevity of notation, we set the following definition:
\begin{df}
Let $(A, \star, \alpha)$ be a hom-associative $k$-algebra. We say that $A$ is a \emph{twisting of an associative algebra} if it arises by the procedure in Ex. \ref{generalized yau-example} from an associative structure $(A, \cdot)$ on the underlying set $A$.
\end{df}
Let for the rest of the section $(A, \star, \alpha, c)$ be a weakly left unital hom-associative algebra with weak left unit $c \in A$, bijective alpha and let $\beta := \alpha^{-1}$. As we already indicated, our aim is to show that $A$ is a twisting of a left unital associative algebra. In order to do this, we show first a few lemmata:
\begin{lem} \label{induced type I3 lemma}
$(A, \star, \beta)$ satisfies the identity $(\beta(x) \star y) \star z = x \star (y \star \beta(z))$.
\end{lem}
\begin{proof}
We calculate
\begin{displaymath}
(\beta(x) \star y) \star z = (\beta(x) \star y) \star \alpha(\beta(z)) = \alpha(\beta(x)) \star (y \star \beta(z)) = x \star (y \star \beta(z))
\end{displaymath}
for all $x,y,z \in A$.
\end{proof}
\begin{lem} \label{symmetry of weak left unit}
The weak left unit obeys the symmetry condition
\begin{displaymath}
(c \star x) \star y=(x \star c) \star y
\end{displaymath}
for all $x,y \in A$.
\end{lem}
\begin{proof}
We have 
\begin{displaymath}
(c \star x) \star y \overset{\alpha \circ \beta = id}{=} (c \star x) \star (c \star \beta(y)) \overset{\text{Hom}}{=} (x \star c) \star (c \star \beta(y)) = (x \star c) \star y.
\end{displaymath}
\end{proof}
Note now that in Example \ref{yau-example} due to Yau, $\alpha$ was supposed an algebra endomorphism for an associative structure $(A, \cdot)$ on $A$. If this is the case, then $\beta$ will also be an algebra endomorphism, i.e. we have $\beta(x \cdot y) = \beta(x) \cdot \beta(y)$ for all $x,y \in A$. In terms of the hom-associative multiplication $\star$, this means that
\begin{displaymath}
\beta(x \star y) = x \cdot y,
\end{displaymath}
and therefore
\begin{displaymath}
\beta(x \star y) = \beta(x) \star \beta(y).
\end{displaymath}
The following lemma provides a generalization of this observation to our situation:
\begin{lem} \label{product rule lemma}
We have
\begin{displaymath}
\beta(x) \star \beta(y) = \beta(c) \star \beta(\beta(x \star y))
\end{displaymath}
for all $x,y \in A$.
\end{lem}
\begin{proof}
We start by noticing that
\begin{equation} \label{eq-a}
x \star \beta(y) = \beta((\beta(c) \star x) \star y)
\end{equation}
for all $x,y \in A$ since
\begin{displaymath}
\beta((\beta(c) \star x) \star y) = \beta(c \star (x \star \beta(y))) = x \star \beta(y)
\end{displaymath}
due to the $(\beta(x) \star y) \star z = x \star (y \star \beta(z))$ for all $x,y,z \in A$ by Lemma \ref{induced type I3 lemma} and $\beta \circ \alpha = id$. Using Eq. (\ref{eq-a}) on $\beta(x) \star \beta(y)$, we get
\begin{equation} \label{eq-b}
\beta(x) \star \beta(y) = \beta((\beta(c) \star \beta(x)) \star y).
\end{equation}
On the other hand, we can also evaluate Eq. (\ref{eq-a}) on $\beta(c)$ and $\beta(\beta(x \star y))$ to obtain
\begin{equation} \label{eq-c}
\beta(c) \star \beta(\beta(x \star y)) = \beta((\beta(c) \star \beta(c)) \star \beta(x \star y)).
\end{equation}
Our goal will now be to prove that the righthand sides of Eq. (\ref{eq-b}) and\\ Eq. (\ref{eq-c}) are the same. By Lemma \ref{symmetry of weak left unit} we see first $$(c \star (\beta(x) \star y)) \star z = ((\beta(x) \star y) \star c) \star z$$ for any $x,y,z \in A$, which using Lemma \ref{induced type I3 lemma} on the righthand side yields
\begin{equation} \label{eq-d}
(c \star (\beta(x) \star y)) \star z = (x \star (y \star \beta(c))) \star z.
\end{equation}
We can use this equation on the argument of $\beta((c \star (\beta(\beta(c)) \star x)) \star y)$ to see
\begin{eqnarray} \label{eq-e}
\beta((c \star (\beta(\beta(c)) \star x)) \star y) &  \overset{\text{Eq. (\ref{eq-d})}}{=} & \beta((\beta(c) \star (x \star \beta(c))) \star y) \nonumber \\ & \overset{\text{Lemma \ref{induced type I3 lemma}}}{=} & \beta(c \star (( x \star \beta(c)) \star \beta(y))) \nonumber \\ &  =  & (x \star \beta(c)) \star \beta(y).
\end{eqnarray}
On the other hand, we also see the more general identity
\begin{eqnarray} \label{eq-f}
\beta(((\beta(\beta(c)) \star x) \star y) \star z) & \overset{\text{Lemma \ref{induced type I3 lemma}}}{=} & \beta((\beta(c) \star (x \star \beta(y))) \star z) \nonumber\\
& \overset{\text{Eq. (\ref{eq-a})}}{=} & (x \star \beta(y)) \star \beta(z). 
\end{eqnarray}
We now show
\begin{equation} \label{eq-g}
(x \star \beta(c)) \star \beta(y \star z) = (x \star \beta(y)) \star z
\end{equation}
by calculating
\begin{eqnarray*}
(x \star \beta(y)) \star z & = & (x \star \beta(y)) \star \beta(c \star z)\\
& \overset{\text{Eq. (\ref{eq-f})}}{=} & \beta(((\beta(\beta(c)) \star x) \star y) \star (c \star z))\\
& \overset{cx = \alpha(x)}{=}  & \beta((c \star (\beta(\beta(c)) \star x)) \star (y \star z))\\ & \overset{\text{Eq. (\ref{eq-e})}}{=} & (x \star \beta(c)) \star \beta(y \star z).
\end{eqnarray*}
Applying now Eq. (\ref{eq-g}) to the righthand side of Eq. (\ref{eq-c}) we see
\begin{eqnarray*}
\beta(c) \star (\beta(\beta(x \star y))) & = & \beta((\beta(c) \star \beta(c)) \star \beta(x \star y))\\
& \overset{\text{Eq. (\ref{eq-g})}}{=} & \beta((\beta(c) \star \beta(x)) \star y)\\
& \overset{\text{Lemma \ref{induced type I3 lemma}}}{=} & \beta(c \star (\beta(x) \star \beta(y))) = \beta(x) \star \beta(y)
\end{eqnarray*}
as desired.
\end{proof}
We are now ready for the proof of the following important lemma, which will be the main ingredient in the proof of the statement we are aiming for. Because of the importance of this result for this section, we re-state explicitly all assumptions about $A$:
\begin{lem} \label{weak unit main tool}
Let $(A, \star, \alpha, c)$ be a weakly left unital hom-associative $k$-algebra with $\alpha$ bijective and let $\beta := \alpha^{-1}$. Then $(A, \star, \beta)$ satisfies the identity $$x \star \beta(y \star z) = \beta(x \star y) \star z$$ for all $x,y,z \in A$.
\end{lem}
\begin{proof}
We see using Lemma \ref{product rule lemma} that
\begin{eqnarray*}
\beta(x \star y) \star z & = & \beta(x \star y) \star \beta(c \star z) \\
& \overset{\text{Lemma \ref{product rule lemma}}}{=} & \beta(c) \star \beta(\beta((x \star y) \star (c \star z))) \\
& \overset{\text{hom-ass}}{=} & \beta(c) \star \beta(\beta((c \star x) \star (y \star z))) \\
& \overset{\text{Lemma \ref{product rule lemma}}}{=} & \beta(c \star x) \star \beta(y \star z) = x \star \beta(y \star z).
\end{eqnarray*}
\end{proof}
We now arrive at the following:
\begin{prop} \label{weakly unital algebras are twistings}
Let $(A, \star, \alpha, c)$ be a left weakly unital hom-associative algebra with bijective $\alpha$. Then there exists a bilinear multiplication $\cdot: A \times A \rightarrow A$ such that $(A, \cdot, c)$ is a left unital associative algebra and such that $(A, \star, \alpha, c)$ arises from $(A, \cdot, c)$ by the procedure of Ex. \ref{generalized yau-example}.
\end{prop}
\begin{proof}
Set $\beta := \alpha^{-1}$ and $x \cdot y := \beta(x \star y)$. It is then clear that $x \star y = \alpha(x \cdot y)$ for all $x,y \in A$ and that $c \cdot x = \beta(c \star x) = x$ for all $x \in A$. So one only needs to show that $(A, \cdot)$ is associative. But this is now easily proven by calculating
\begin{displaymath}
x \cdot (y \cdot z) = \beta(x \star \beta(y \star z)) \overset{\text{Lemma \ref{weak unit main tool}}}{=} \beta(\beta(x \star y) \star z) = (x \cdot y) \cdot z.
\end{displaymath}
\end{proof}

 \noindent {\bf Acknowledgements.}  The authors were supported in their work by research grant R1F105L15 of the University of Luxembourg (Martin Schlichenmaier).\\
This work was aided by use of the first order logic theorem prover Prover9 and the finite countermodel searcher Mace4. The most important contribution to be mentioned in this context is without doubt the proof of Lemma \ref{weak unit main tool} including Lemma \ref{product rule lemma}, which is essentially a structured translation to natural language of a proof discovered by Prover9.\\
The first author also wants to thank the Max Planck Institute for Mathematics in Bonn for excellent working conditions.\\
We also thank Abdenacer Makhlouf and Donald Yau for useful comments on previous work in this direction.

\vskip 1cm



\begin{thebibliography}{99}



\bibitem{AMS} H. Ataguema, A. Makhlouf, S. Silvestrov, \emph{Generalization of n-ary Nambu algebras and
beyond} (Preprint: arXiv:0812.4058).
\bibitem{CG09} S. Caenepeel, I. Goyvaerts, \emph{Hom-Hopf algebras} arXiv:0907.0187v1.
\bibitem{GR04} M. Goze, E. Remm,
\emph{Lie-admissible algebras and operads}, J. Algebra \textbf{273}
(2004), 129--152.
\bibitem{HLS}  J. T. Hartwig, D. Larsson, S.D. Silvestrov,
\emph{Deformations of Lie algebras using $\sigma$-derivations}, J.
Algebra \textbf{295} (2006), 314-361.
\bibitem{LS1} D. Larsson, S. D. Silvestrov,
\emph{Quasi-hom-Lie algebras, Central Extensions and 2-cocycle-like
identities}, J. Algebra \textbf{288} (2005), 321--344.
\bibitem{LS2} D. Larsson, S. D. Silvestrov, \emph{Quasi-Lie algebras}, in "Noncommutative
Geometry and Representation Theory in Mathematical Physics",
Contemp. Math., 391, Amer. Math. Soc., Providence, RI, (2005),
241-248.
\bibitem{LS3} D. Larsson, S. D. Silvestrov,
\emph{Quasi-deformations of $sl_2(\mathbb{F})$ using twisted
derivations}, to appear in Communications in Algebra. (Preprint:
arXiv:math/0506172).
\bibitem{MSI} A. Makhlouf, S. D. Silvestrov,
\emph{Hom-algebra structures},  J. Gen. Lie Theory Appl. Vol
\textbf{2} (2) 2008, pp 51-64.
\bibitem{MSII} A. Makhlouf, S. D. Silvestrov, \emph{Hom-Algebras and Hom-Coalgebras}, Preprints in
Mathematical Sciences, Lund University, Centre for Mathematical
Sciences, Centrum Scientiarum Mathematicarum, (2008) ;
arXiv:0811.0400 [math.RA] (2007).
\bibitem{MSIII} A. Makhlouf, S. D. Silvestrov, \emph{Hom-Lie admissible Hom-coalgebras
and Hom-Hopf algebras}, Published as Chapter 17, pp 189-206, S.
Silvestrov, E. Paal, V. Abramov, A. Stolin, (Eds.), Generalized Lie
theory in Mathematics, Physics and Beyond, Springer-Verlag, Berlin,
Heidelberg, (2008).
\bibitem{MSIV} A. Makhlouf, S. D. Silvestrov, \emph{Notes on
Formal deformations of Hom-Associative and Hom-Lie algebras}, To
appear in Forum Mathematicum.
\bibitem{Yau:EnvLieAlg} D. Yau,
\emph{Enveloping algebra of Hom-Lie algebras}, J. Gen. Lie Theory
Appl. Vol \textbf{2} (2) 2008, pp 95-108.
\bibitem{Yau:HOMology} D. Yau,
\emph{Hom-algebras and homology}, arXiv:0712.3515v2
\bibitem{Yau:Hom-bi} D. Yau,\emph{Hom-bialgebras and comodule algebras}, arXiv:0810.4866v1
[math.RA] (2008).
\bibitem{Yau:Hom-Yang} D. Yau, \emph{The Hom-Yang-Baxter equation, Hom-Lie algebras, and quasi-triangular bialgebras}
Journal of Physics A 42 (2009).
\bibitem{Yau:Hom-qgI} D. Yau, \emph{Hom-quantum groups II: quasi-triangular Hom-bialgebras}, arXiv:0906.4128v1 [math-ph] (2009).
\bibitem{Yau:Hom-qgII} D. Yau, \emph{Hom-quantum groups II: cobraided Hom-bialgebras and Hom-quantum geometry}, arXiv:0907.1880v1 [math.QA] (2009).
\end{thebibliography}
\end{document}